\documentclass[a4j,12pt,draft]{article}
\usepackage{indent first}
\usepackage[a4paper,vmargin={1in,1in},hmargin={2.4cm,2.4cm}]{geometry}
\usepackage{amsmath}
\usepackage{amssymb}
\usepackage{amsthm}
\usepackage{mathrsfs}
\usepackage{titlesec}
\titleformat{\section}{\large\bfseries}{\thesection}{1em}{}
\numberwithin{equation}{section}

\theoremstyle{definition}
\newtheorem{definition}{Definition}[section]
\theoremstyle{plain}
\newtheorem{thm}{Theorem}
\theoremstyle{plain}
\newtheorem{prop}{Proposition}[section]
\theoremstyle{plain}
\newtheorem{lem}[prop]{Lemma}
\theoremstyle{plain}
\newtheorem{cor}[thm]{Corollary}
\theoremstyle{definition}
\newtheorem{rmk}[prop]{Remark}

\begin{document}

\begin{center}
\Large{\textbf{On the Zeros of the Second Derivative of the Riemann Zeta Function under the Riemann Hypothesis}}
\end{center}

\vspace{2mm}
\begin{center}
\large{Ade Irma Suriajaya}
\let\thefootnote\relax\footnote{\emph{Date}: September 27, 2013.}
\let\thefootnote\relax\footnote{The author is supported by Nitori International Scholarship Foundation.}
\end{center}

\begin{center}
Graduate School of Mathematics, Nagoya University,\\
Furo-cho, Chikusa-ku, Nagoya, 464-8602, Japan
\end{center}

\vspace{3mm}
\begin{abstract}
The number of zeros and the distribution of the real part of non-real zeros of the derivatives of the Riemann zeta function have been investigated by Berndt, Levinson, Montgomery, and Akatsuka. Berndt, Levinson, and Montgomery investigated the general case, meanwhile Akatsuka gave sharper estimates for the first derivative of the Riemann zeta function under the truth of the Riemann hypothesis. In this paper, we extend the results of Akatsuka to the second derivative of the Riemann zeta function.
\end{abstract}

\vspace{1mm}
\section{Introduction}
\label{sec:1}

The theory of the Riemann zeta function $\zeta(s)$ has been studied for over 150 years. Among the topics of research, the study of its zeros has been one of the main focus of research. Lately, the study of the zeros of its derivatives has also been part of the research area.
In fact, in 1970, Berndt \cite[Theorem]{ber} proved that
\begin{equation} \label{eq:nkt}
N_k(T) = \frac{T}{2\pi}\log{\frac{T}{4\pi}}-\frac{T}{2\pi}+O(\log{T})
\end{equation}
where $N_k(T)$ denotes the number of zeros of the $k$-th derivative of the Riemann zeta function, we write as $\zeta^{(k)}(s)$, with $0<\text{Im}\,(s)\leq T$, counted with multiplicity, for any positive integer $k\geq1$.
And in 1974, Levinson and Montgomery \cite[Theorem 10]{lev} showed that for any positive integer $k\geq1$,
\begin{equation} \label{eq:kdist}
\begin{aligned}
\sum_{\substack{\rho^{(k)}=\beta^{(k)}+i\gamma^{(k)},\\ \zeta^{(k)}(\rho^{(k)})=0,\, 0<\gamma^{(k)}\leq T}} \left(\beta^{(k)}-\frac{1}{2}\right) = \frac{kT}{2\pi}\log{\log{\frac{T}{2\pi}}} &+ \frac{1}{2\pi}\left(\frac{1}{2}\log{2} - k\log{\log{2}}\right)T\\
&\quad\quad\quad\quad\quad\quad\quad- k\text{Li}\left(\frac{T}{2\pi}\right) + O(\log{T})
\end{aligned}
\end{equation}
where the sum is counted with multiplicity and
\begin{equation*}
\text{Li}(x) := \int_2^x \frac{dt}{\log{t}}.
\end{equation*}
In 2012, Akatsuka \cite[Theorem 1 and Theorem 3]{aka} improved each of the error term of equations \eqref{eq:nkt} and \eqref{eq:kdist} for the case of $k=1$ under the assumption of the truth of the Riemann hypothesis.
His results are
\begin{equation*}
\begin{aligned}
\sum_{\substack{\rho'=\beta'+i\gamma',\\ \zeta'(\rho')=0,\, 0<\gamma'\leq T}} \left(\beta'-\frac{1}{2}\right) = \frac{T}{2\pi}\log{\log{\frac{T}{2\pi}}} &+ \frac{1}{2\pi}\left(\frac{1}{2}\log{2} - \log{\log{2}}\right)T\\
&\quad\quad\quad\quad\quad\quad\quad- \text{Li}\left(\frac{T}{2\pi}\right) + O((\log{\log{T})^2})
\end{aligned}
\end{equation*}
and
\begin{equation*}
N_1(T) = \frac{T}{2\pi}\log{\frac{T}{4\pi}}-\frac{T}{2\pi}+O\left(\frac{\log{T}}{(\log{\log{T}})^{1/2}}\right)
\end{equation*}
if the Riemann hypothesis is true.
Here, we extend these results to the case of $k=2$.
Before we introduce our results, we define some notations which are going to be used throughout this paper.

\vspace{3mm}
In this paper we denote by $\mathbb{Z}$, $\mathbb{R}$, and $\mathbb{C}$
the set of all rational integers, the set of all real numbers, and the set of all complex numbers, respectively.
Let $\rho = \beta + i\gamma$ and $\rho'' = \beta'' + i\gamma''$ represent the nontrivial zeros of the Riemann zeta function and the non-real zeros of the second derivative of the Riemann zeta function, respectively. Then we define $N(T)$ and $N_2(T)$ as follows.

\begin{definition}
For $T>0$, we define
\begin{equation*}
N(T) := \sharp'\{\rho=\beta+i\gamma \,|\, 0<\gamma\leq T\}
\end{equation*}
and
\begin{equation*}
N_2(T) := \sharp'\{\rho''=\beta''+i\gamma'' \,|\, 0<\gamma''\leq T\}
\end{equation*}
where $\sharp'$ means the number of elements counted with multiplicity.
\end{definition}
Next, for any complex number $s$, we write Re$(s)$ and Im$(s)$ as $\sigma$ and $t$, respectively.
Finally, we abbreviate the Riemann Hypothesis as RH.

Below we state our results, each of which is analogous to Theorem 1, Corollary 2, and Theorem 3 of \cite{aka}, respectively.

\begin{thm} \label{cha1}
Assume RH. Then for any $T>2\pi$, we have
\begin{equation*}
\begin{aligned}
\sum_{\substack{\rho''=\beta''+i\gamma'',\\ 0<\gamma''\leq T}} \left(\beta''-\frac{1}{2}\right) = \frac{2T}{2\pi}\log{\log{\frac{T}{2\pi}}} &+ \frac{1}{2\pi}\left(\frac{1}{2}\log{2} - 2\log{\log{2}}\right)T\\
&\quad\quad\quad\quad\quad\quad\quad- 2\emph{Li}\left(\frac{T}{2\pi}\right) + O((\log{\log{T}})^2).
\end{aligned}
\end{equation*}
\end{thm}

\begin{cor} \label{cha12} (Cf. \cite[Theorem 3]{lev}.)
Assume RH. Then for $0<U<T$ (where $T$ is restricted to satisfy $T>2\pi$), we have
\begin{equation*}
\begin{aligned}
\sum_{\substack{\rho''=\beta''+i\gamma'',\\ T<\gamma''\leq T+U}} \left(\beta''-\frac{1}{2}\right) = \frac{2U}{2\pi}\log{\log{\frac{T}{2\pi}}} &+ \frac{1}{2\pi}\left(\frac{1}{2}\log{2} - 2\log{\log{2}}\right)U\\
&\quad\quad\quad\quad\quad\quad\quad+ O\left(\frac{U^2}{T\log{T}}\right) + O((\log{\log{T}})^2).
\end{aligned}
\end{equation*}
\end{cor}

\begin{thm} \label{cha2}
Assume RH. Then for $T\geq2$, we have
\begin{equation*}
N_2(T) = \frac{T}{2\pi}\log{\frac{T}{4\pi}} - \frac{T}{2\pi} + O\left(\frac{\log{T}}{(\log{\log{T}})^{1/2}}\right).
\end{equation*}
\end{thm}

\vspace{2mm}
Before we begin the next section, we intend to give a brief outline of the proofs.
Nevertheless, since the steps of our proofs basically follow those given in \cite{aka} with a few crucial modifications, instead of the outline of the proof, we present the main modifications related to the proofs.
We note here that Lemma 2.1 of \cite{aka} plays a really important role in the proofs and so most of the points which are modified are parts of that lemma.

First of all, condition 2 of the lemma is related to the functional equation for $\zeta'(s)$. Hence in our case, we need to consider for $\zeta''(s)$ which gives us more terms to estimate, some of which are not logarithmic derivatives of some functions that need to be taken care differently from those in the case of $\zeta'(s)$.

Next is condition 4 of the lemma. For $\zeta'(s)$, the term we need to estimate was $\frac{\zeta'}{\zeta}(s)$ which is just the logarithmic derivative of $\zeta(s)$. In \cite{aka}, the inequality $\text{Re}\,\left(\frac{\zeta'}{\zeta}(s)\right)<0$ was obtained, however for $\zeta''(s)$, the sign of $\text{Re}\,\left(\frac{\zeta''}{\zeta}(s)\right)$ does not seem to stay unchanged in any region defined by $x\leq\sigma<\frac{1}{2},\, t\geq y$ for some $x\leq-1$ and large $y>0$. Since it is sufficient to show that $\frac{\zeta''}{\zeta}(s)$ is holomorphic and non-zero, and has bounded argument in some such region, we modify the condition in such a way.

Finally, the region $\frac{1}{2}<\sigma\leq a$ considered in Lemma 2.3 of \cite{aka} does not work well for $\frac{\zeta''}{\zeta}(s)$. The reason is that the current best estimation of $\frac{\zeta''}{\zeta}(s)$ depends on the usage of Cauchy's integral formula, hence we need to keep a certain distance between $\frac{1}{2}$ and the infimum of $\sigma$ in the region. Therefore, we put here a small distance $\epsilon_0>0$.


\section{Proof of Theorem \ref{cha1} and Corollary \ref{cha12}}
\label{sec:2}

In this section we give the proof of Theorem 1 and Corollary 2. For that purpose, we need a few lemmas and a proposition which are analogues of those in \cite{aka}.
For convenience, we define two functions $F(s)$ and $G_2(s)$ as follows.

\begin{definition}
\begin{equation*}
F(s) := 2^s\pi^{s-1}\sin{\left(\frac{\pi s}{2}\right)}\Gamma(1-s),\quad G_2(s) := \frac{2^s}{(\log{2})^2}\zeta''(s).
\end{equation*}
\end{definition}

\vspace{3mm}
By the above definition of $F(s)$, we can check easily that the functional equation for $\zeta(s)$ states
\begin{equation} \label{eq:fe}
\zeta(s) = F(s)\zeta(1-s).
\end{equation}

To begin, we introduce the following lemma which is the analogue of Lemma 2.1 of \cite{aka} for the case of $\zeta''(s)$.

\vspace{3mm}
\begin{lem} \label{lem1}
Assume RH. Then there exist $10\leq a\leq14$, $\sigma_0\leq-1$, and $t_0\geq\max\{30,-\sigma_0\}$ such that
\begin{equation*}
\begin{aligned}
&1. \left|G_2(s)-1\right| < \frac{1}{2}\left(\frac{2}{3}\right)^{\sigma/2}, \quad \text{for any}\,\,\sigma\geq a;\\
&2. \left|2\frac{1}{\frac{F''}{F'}(s)}\frac{\zeta'}{\zeta}(1-s) - \frac{1}{\frac{F''}{F}(s)}\frac{\zeta''}{\zeta}(1-s)\right| \leq 2^\sigma, \quad \text{for}\,\, \sigma\leq\sigma_0\,\, \text{and}\,\, t\geq2;\\
&3. \left|\frac{F''}{F}(s)\right|\geq1\,\, \text{and there exists an}\,\, n\in\mathbb{Z} \,\text{such that}\, \arg{\frac{F''}{F}(s)}\in \left[-\frac{\pi}{6} + 2n\pi,\frac{\pi}{6} + 2n\pi\right],\\
&\quad\quad\quad\quad\quad\quad\quad\quad\quad\quad\quad\quad\quad\quad\quad\quad\quad\quad\quad\quad\quad\quad\quad\quad \text{for}\,\, \sigma_0\leq\sigma\leq\frac{1}{2}\,\, \text{and}\,\, t\geq t_0-1;\\
&4.\, \frac{\zeta''}{\zeta}(s)\neq0,\,\, \text{and there exists an}\,\, n\in\mathbb{Z} \,\text{such that}\, \arg{\frac{\zeta''}{\zeta}(s)}\in \left[(2n-1)\pi,(2n+1)\pi\right],\\
&\quad\quad\quad\quad\quad\quad\quad\quad\quad\quad\quad\quad\quad\quad\quad\quad\quad\quad\quad\quad\quad\quad\quad\quad \text{for}\,\, \sigma_0\leq\sigma<\frac{1}{2}\,\, \text{and}\,\, t\geq t_0-1;\\
&5.\, \zeta(\sigma+it_0)\neq0,\, \zeta''(\sigma+it_0)\neq0, \quad \text{for any}\,\, \sigma\in\mathbb{R}.\\
\end{aligned}
\end{equation*}
\end{lem}

\vspace{1mm}
\begin{proof}
\raggedright
$\quad$\\

\vspace{2mm}
1. Let $a=12$, we show that this satisfies condition 1.
For $\sigma\geq a\,\, (=12)$, we have
\begin{equation*}
G_2(s) = \frac{2^s}{(\log{2})^2}\zeta''(s) = \frac{2^s}{(\log{2})^2}\sum_{n=1}^\infty \frac{(\log{n})^2}{n^s} = 1 + \left(\frac{\log{3}}{\log{2}}\right)^2 \left(\frac{2}{3}\right)^s + \frac{2^s}{(\log{2})^2}\sum_{n=4}^\infty \frac{(\log{n})^2}{n^s}.
\end{equation*}
Hence,
\begin{equation*}
\begin{aligned}
&\left|G_2(s) - 1\right| \leq \left|\left(\frac{\log{3}}{\log{2}}\right)^2 \left(\frac{2}{3}\right)^s\right| + \left|\frac{2^s}{(\log{2})^2}\sum_{n=4}^\infty \frac{(\log{n})^2}{n^s}\right|\\
&\,\,\,\leq \left(\frac{\log{3}}{\log{2}}\right)^2 \left(\frac{2}{3}\right)^\sigma + \frac{2^\sigma}{(\log{2})^2}\int_3^\infty \frac{(\log{x})^2}{x^\sigma}dx\\
&\,\,\,\leq \left\{\left(\frac{\log{3}}{\log{2}}\right)^2 + \frac{3}{\sigma-1}\left(\frac{\log{3}}{\log{2}}\right)^2 + \frac{6\log{3}}{(\sigma-1)^2(\log{2})^2} + \frac{6}{(\sigma-1)^3(\log{2})^2}\right\} \left(\frac{2}{3}\right)^\sigma\\
&\,\,\,\leq \left\{\left(\frac{\log{3}}{\log{2}}\right)^2 + \frac{3}{a-1}\left(\frac{\log{3}}{\log{2}}\right)^2 + \frac{6\log{3}}{(a-1)^2(\log{2})^2} + \frac{6}{(a-1)^3(\log{2})^2}\right\} \left(\frac{2}{3}\right)^{a/2} \left(\frac{2}{3}\right)^{\sigma/2}\\
&\,\,\,\leq \left\{\left(\frac{\log{3}}{\log{2}}\right)^2 + \frac{3}{11}\left(\frac{\log{3}}{\log{2}}\right)^2 + \frac{6\log{3}}{11^2(\log{2})^2} + \frac{6}{11^3(\log{2})^2}\right\} \left(\frac{2}{3}\right)^6 \left(\frac{2}{3}\right)^{\sigma/2}< \frac{1}{2} \left(\frac{2}{3}\right)^{\sigma/2}.
\end{aligned}
\end{equation*}
Thus, we have proven the existence of such constant $a$.\\

\vspace{2mm}
2. Here we consider for $\sigma\leq-1$ and $t\geq2$.

Firstly we estimate $\frac{F''}{F}(s)$.
We apply Stirling's formula to $\frac{\Gamma'}{\Gamma}(1-s)$ in the region $-\pi+\delta\leq\arg{(1-s)}\leq\pi-\delta$ by taking $\delta$ to be $\frac{\pi}{2}$ which is sufficient for the region we are considering. Then we have
\begin{equation} \label{eq:flog}
\begin{aligned}
\frac{F'}{F}(s) &= \log{2\pi} + \frac{\pi}{2}\cot{\left(\frac{\pi s}{2}\right)} - \frac{\Gamma'}{\Gamma}(1-s)\\
&= \log{2\pi} + \frac{\pi}{2}\cot{\left(\frac{\pi s}{2}\right)} - \log{(1-s)} + \frac{1}{2(1-s)} + \int_0^\infty \frac{[u]-u+\frac{1}{2}}{(u+1-s)^2} du.
\end{aligned}
\end{equation}
Now apply the product rule of derivatives to obtain
\begin{equation} \label{eq:f2}
\frac{F''}{F}(s) = \left(\frac{F'}{F}(s)\right)' + \left(\frac{F'}{F}(s)\right)^2.
\end{equation}
Using equation \eqref{eq:flog} and noting that $\cot{\left(\frac{\pi s}{2}\right)} = O(1)$ for $t\geq2$, we can show that
\begin{equation*}
\left|\frac{F''}{F}(s)\right| \geq (\log{|1-s|})^2 - | O(\log{|1-s|})|.
\end{equation*}
Thus, we can take $\sigma_1\leq-1$ sufficiently small (i.e. sufficiently large in the negative direction) so that $|1-s|$ is large enough such that
\begin{equation} \label{eq:21}
\left|\frac{F''}{F}(s)\right| \geq \frac{1}{4}(\log{|1-s|})^2 \quad (^\forall\sigma\leq\sigma_1,\,\, t\geq2).
\end{equation}

Next, we estimate $\frac{F''}{F'}(s)$.
We note that this is the logarithmic derivative of $F'(s)$, then apply Stirling's formula to $\frac{\Gamma'}{\Gamma}(1-s)$ as in equation \eqref{eq:flog} and by noting that $\csc{\left(\frac{\pi s}{2}\right)} = O(1),\, \cot{\left(\frac{\pi s}{2}\right)} = O(1)$ when $t\geq2$, we can show that
\begin{equation*}
\frac{F''}{F'}(s) = -\log{(1-s)} + O(1)
\end{equation*}
for small $\sigma\leq-1$.
Hence we can take a $\sigma_2\leq-1$ such that
\begin{equation} \label{eq:22}
\left|\frac{F''}{F'}(s)\right| \geq \frac{1}{2}\log{(1-\sigma)} \quad (\sigma\leq\sigma_2,\, t\geq2).
\end{equation}

Now we estimate $\frac{\zeta'}{\zeta}(1-s)$. For $\sigma\leq-1$, we have
\begin{equation} \label{eq:23}
\begin{aligned}
\left|\frac{\zeta'}{\zeta}(1-s)\right|
&= \left| -\sum_{n=1}^\infty \frac{\Lambda(n)}{n^{1-s}}\right|
\leq \sum_{n=1}^\infty \frac{\Lambda(n)}{n^{1-\sigma}}
= \frac{\log{2}}{2^{1-\sigma}} + \sum_{n=3}^\infty \frac{\Lambda(n)}{n^{1-\sigma}}
\leq \frac{\log{2}}{2^{1-\sigma}} + \sum_{n=3}^\infty \frac{\log{n}}{n^{1-\sigma}}\\
&\leq \frac{\log{2}}{2^{1-\sigma}} + \int_2^\infty \frac{\log{x}}{x^{1-\sigma}} dx
\leq 2^\sigma\left( \frac{3}{2}\log{2}+1 \right).
\end{aligned}
\end{equation}

Finally, we estimate $\frac{\zeta''}{\zeta}(1-s)$. For $\sigma\leq -1$, by also using the result of equation \eqref{eq:23} we can calculate
\begin{equation} \label{eq:24}
\begin{aligned}
\left|\frac{\zeta''}{\zeta}(1-s)\right| &= \left| -\left(\frac{\zeta'}{\zeta}(1-s)\right)' + \left(\frac{\zeta'}{\zeta}(1-s)\right)^2 \right|
=\left| \left(\sum_{n=1}^\infty \frac{\Lambda(n)}{n^{1-s}}\right)' +  \left(\sum_{n=1}^\infty \frac{\Lambda(n)}{n^{1-s}}\right)^2 \right| \\
&= \left| \sum_{n=1}^\infty \frac{\Lambda(n)\log{n}}{n^{1-s}} +  \left(\sum_{n=1}^\infty \frac{\Lambda(n)}{n^{1-s}}\right)^2 \right|
\leq \sum_{n=1}^\infty \frac{\Lambda(n)\log{n}}{n^{1-\sigma}} +  \sum_{n=1}^\infty \frac{\Lambda(n)}{n^{1-\sigma}} \sum_{n=1}^\infty \frac{\Lambda(n)}{n^{1-\sigma}} \\
&\leq 2^\sigma\left( \frac{19}{8}(\log{2})^2 + \frac{13}{4}\log{2} + \frac{5}{2} \right).
\end{aligned}
\end{equation}

Thus, combining equations \eqref{eq:21} to \eqref{eq:24} and using simple calculations, we can show that for $\sigma\leq\min\{\sigma_1, \sigma_2\}$, we have
\begin{equation*}
\begin{aligned}
\left| 2\frac{1}{\frac{F''}{F'}(s)}\frac{\zeta'}{\zeta}(1-s) - \frac{1}{\frac{F''}{F}(s)}\frac{\zeta''}{\zeta}(1-s)\right|
&\leq 2\left|\frac{1}{\frac{F''}{F'}(s)}\right| \left|\frac{\zeta'}{\zeta}(1-s)\right| + \left|\frac{1}{\frac{F''}{F}(s)}\right| \left|\frac{\zeta''}{\zeta}(1-s)\right| \\
&\leq 32 \frac{2^\sigma}{\log{(1-\sigma)}}.
\end{aligned}
\end{equation*}

Hence, by taking $\sigma_0\leq\min\{\sigma_1, \sigma_2\}$ such that $\log{(1-\sigma)}\geq32$, we obtain
\begin{equation*}
\left| 2\frac{1}{\frac{F''}{F'}(s)}\frac{\zeta'}{\zeta}(1-s) - \frac{1}{\frac{F''}{F}(s)}\frac{\zeta''}{\zeta}(1-s)\right| \leq 2^\sigma, \quad (\text{for any}\,\, \sigma\leq\sigma_0,\, t\geq2).
\end{equation*}

\vspace{2mm}
3. Now with the above $\sigma_0$, we are going to look for $t_0$ for which conditions 3 to 5 hold.
We start by examining condition 3. We consider for $\sigma_0\leq\sigma\leq\frac{1}{2}$ and $t\geq29$.
As in the proof of condition 2, we apply Stirling's formula in the region $-\frac{\pi}{2}\leq\arg{(1-s)}\leq\frac{\pi}{2}$ which is sufficient for our purpose to obtain equation \eqref{eq:flog}. Thus, by again using the product rule of derivatives for $\frac{F'}{F}(s)$ (eq. \eqref{eq:f2}), for $\sigma_0\leq\sigma\leq\frac{1}{2}$ and $t\geq29$ we have
\begin{equation*}
\frac{F''}{F}(s) = (\log{(1-s)})^2 + O_{\sigma_0}(\log{t}).
\end{equation*}
Consequently,
\begin{equation} \label{eq:f2log}
\left|\frac{F''}{F}(s)\right| \geq |\log{(1-s)}|^2 - |O_{\sigma_0}(\log{t})| \geq (\log{t})^2 - |O_{\sigma_0}(\log{t})|.
\end{equation}
Thus, we can take $t_1\geq30$ for which $(\log{t})^2 - |O_{\sigma_0}(\log{t})| \geq 1$ holds for $^\forall t\geq t_1-1$.
We note from equation \eqref{eq:f2log} that $\lim_{t\rightarrow\infty} \frac{F''}{F}(s) = \infty$ holds uniformly for $\sigma_0\leq\sigma\leq\frac{1}{2}$.
Hence there exists a $t_2\geq30$ such that
\begin{equation*}
-\frac{\pi}{6}+2n\pi \leq \arg{\frac{F''}{F}(s)} \leq \frac{\pi}{6}+2n\pi
\end{equation*}
for some $n\in\mathbb{Z}$ (depends on our choice of logarithmic branch) and for any $s=\sigma+it$ satisfying $\sigma_0\leq\sigma\leq\frac{1}{2}$ and $t\geq t_2-1$.

By the above calculations, we find that $\max\{t_1, t_2, -\sigma_0\}$ is a candidate for $t_0$, thus we have proven that $t_0\geq\max\{30,-\sigma_0\}$ such that condition 3 holds exists. Since we want $t_0$ to also satisfy conditions 4 and 5, we need to examine those conditions to completely prove the existence of $t_0$.

\vspace{2mm}
4. According to \cite[Theorem 1]{yil}, RH implies that $\zeta''(s)$ has no zeros in the strip $0\leq\sigma<\frac{1}{2}$. Furthermore, by \cite[Theorem 2]{yil}, $\zeta''(s)$ has exactly only one pair of non-real zeros in the left half-plane $\sigma<0$. Let us denote that pair of zeros as $\rho''_0$ and $\overline{\rho''_0}$ (the complex conjugate of $\rho''_0$).
Now we take $t_3= \max\{30, |\text{Im}(\rho_0'')|+2 \}$, then $\zeta''(s)\neq0$ in the region defined by $\sigma_0\leq\sigma<\frac{1}{2}$, $t\geq t_3-1$.
Assuming RH, we note that $\zeta(s)$ is holomorphic and has no zeros in that region, hence $\frac{\zeta''}{\zeta}(s)\neq0$ there.

Next we show that there exists an $n\in\mathbb{Z}$ such that $(2n-1)\pi \leq\arg{\frac{\zeta''}{\zeta}(s)}\leq (2n+1)\pi$ holds for $\sigma_0\leq\sigma<\frac{1}{2}$ and $t\geq t_4-1$ for some $t_4\geq30$.
For that purpose, we show that there exists a $t_4\geq30$ such that
\begin{equation*}
\text{Re}\left(\frac{\zeta'}{\zeta}(s)\right) < 0 \quad\text{and}\quad \text{Re}\left(\frac{\zeta''}{\zeta'}(s)\right) < 0
\end{equation*}
hold for $\sigma_0\leq\sigma<\frac{1}{2}$ and $t\geq t_4-1$.
We refer to \cite[pp. 64--65]{lev} and modify some calculations there to show the above claim.
We note that the trivial (real) zeros of $\zeta(s)$ are $s=-2m$ ($m=1,2,3,...$).
Also since we are assuming RH, referring to \cite{spe}, $\zeta'(s)\neq0$ when $\sigma<\frac{1}{2}$ and $t\neq0$. Thus, according to \cite{spi1}, we may denote by $-a_m$ the real zeros of $\zeta'(s)$ where $-a_m\in(-1-2m,1-2m)$ ($m=1,2,3,...$).

Applying the Weierstrass factorization theorem and by RH, we have the well known formulas
\begin{equation}
\begin{aligned}
\frac{\zeta'}{\zeta}(s) &= \log{2\pi} - 1 - \frac{1}{s-1} + \sum_{m=1}^\infty\left(\frac{1}{s+2m}-\frac{1}{2m}\right) + \sum_\gamma\left(\frac{1}{s-\frac{1}{2}-i\gamma}+\frac{1}{\frac{1}{2}+i\gamma}\right),\\
\frac{\zeta''}{\zeta'}(s) &= \frac{\zeta''}{\zeta'}(0) - 2 - \frac{2}{s-1} + \sum_{m=1}^\infty\left(\frac{1}{s+a_m}-\frac{1}{a_m}\right) + \sum_{\rho'}\left(\frac{1}{s-\rho'}+\frac{1}{\rho'}\right).
\end{aligned}
\end{equation}

Thus, for $\sigma_0\leq\sigma<\frac{1}{2}$ and $t\geq29$,
\begin{equation*}
\begin{aligned}
\text{Re}\left(\frac{\zeta'}{\zeta}(s)\right) &= \log{2\pi} - 1 - \frac{\sigma-1}{(\sigma-1)^2+t^2} + \sum_{m=1}^\infty\left(\frac{\sigma+2m}{(\sigma+2m)^2+t^2}-\frac{1}{2m}\right) \\
&\quad\quad\quad\quad\quad\quad\quad\quad\quad\quad\quad\quad\,\,\,\,+ \sum_\gamma\frac{\sigma-\frac{1}{2}}{\left(\sigma-\frac{1}{2}\right)^2+(t-\gamma)^2} + \sum_\gamma\frac{1}{\frac{1}{2}+i\gamma}\\
&\leq \log{2\pi} - 1 - \frac{\sigma-1}{(\sigma-1)^2+t^2} + \sum_{m=1}^\infty\left(\frac{\sigma+2m}{(\sigma+2m)^2+t^2}-\frac{1}{2m}\right) + \sum_\gamma\frac{1}{\frac{1}{2}+i\gamma}\\
&= \sum_{m=1}^\infty\left(\frac{\sigma+2m}{(\sigma+2m)^2+t^2}-\frac{1}{2m}\right) + O_{\sigma_0}(1)\\
&= -|s|^2 \sum_{m=1}^\infty\frac{1}{2m((\sigma+2m)^2+t^2)} - \sigma\sum_{m=1}^\infty\frac{1}{(\sigma+2m)^2+t^2} + O_{\sigma_0}(1) \\
&= -|s|^2 \sum_{m=1}^\infty\frac{1}{2m((\sigma+2m)^2+t^2)} + O_{\sigma_0}(1)
\leq -\frac{2}{9}\log{|s|} + O_{\sigma_0}(1).
\end{aligned}
\end{equation*}
Therefore, we can take $t'_4\geq30$ for which $|O_{\sigma_0}(1)|<\frac{2}{9}\log{|s|}$ holds for $\sigma_0\leq\sigma<\frac{1}{2}$ and $t\geq t'_4-1$.
Consequently, $\text{Re}\left(\frac{\zeta'}{\zeta}(s)\right)<0$ for $\sigma_0\leq\sigma<\frac{1}{2}$ and $t\geq t'_4-1$.

Similarly for $\text{Re}\left(\frac{\zeta''}{\zeta'}(s)\right)$, using the fact that $-a_m\in(-1-2m,1-2m)$ for $m=1,2,3,...$ and that $\frac{1}{2}\leq\beta'<3$ (cf. \cite{spe} and \cite[p. 678]{spi2}), we can also show that
\begin{equation*}
\text{Re}\left(\frac{\zeta''}{\zeta'}(s)\right) \leq -\frac{2}{9}\log{|s|} + O_{\sigma_0}(1).
\end{equation*}
Therefore, we can also take $t''_4\geq30$ for which $\text{Re}\left(\frac{\zeta''}{\zeta'}(s)\right)<0$ holds for $\sigma_0\leq\sigma<\frac{1}{2}$ and $t\geq t''_4-1$.

We set $t_4:=\max\{t'_4,t''_4\}$, then
$\text{Re}\left(\frac{\zeta'}{\zeta}(s)\right)<0$ and $\text{Re}\left(\frac{\zeta''}{\zeta'}(s)\right)<0$ hold for $\sigma_0\leq\sigma<\frac{1}{2}$ and $t\geq t_4-1$.
Hence, there exist $k,\, l\in\mathbb{Z}$ for which
\begin{equation*}
\frac{\pi}{2}+2k\pi\leq\arg{\frac{\zeta'}{\zeta}(s)}\leq\frac{3\pi}{2}+2k\pi,\quad \frac{\pi}{2}+2l\pi\leq\arg{\frac{\zeta''}{\zeta'}(s)}\leq\frac{3\pi}{2}+2l\pi
\end{equation*}
hold for $\sigma_0\leq\sigma<\frac{1}{2},\, t\geq t_4-1$.

Consequently,
\begin{equation*}
(2(k+l)+1)\pi \leq \arg{\frac{\zeta''}{\zeta}(s)} \leq (2(k+l)+3)\pi, \quad (\sigma_0\leq\sigma<\frac{1}{2},\, t\geq t_4-1).
\end{equation*}
Letting $n:=k+l+1$, we obtain
\begin{equation*}
(2n-1)\pi \leq \arg{\frac{\zeta''}{\zeta}(s)} \leq (2n+1)\pi, \quad (\sigma_0\leq\sigma<\frac{1}{2},\, t\geq t_4-1).
\end{equation*}

Combining the proof of condition 3 and the above calculations, we find that $\max\{t_1, t_2, t_3, t_4, -\sigma_0\}$ is a candidate for $t_0$, thus we have proven that $t_0\geq\max\{30,-\sigma_0\}$ for which condition 3 and 4 hold exists.

\vspace{2mm}
5. Now we set $t_5:=\max\{t_1, t_2, t_3, t_4, -\sigma_0\}$.
Referring to \cite[p. 678]{spi2}, we see that $\zeta''(s)\neq0 \quad (\sigma\geq5)$.
Also since $t_5\geq|\text{Im}(\rho_0'')|+2$, as stated in the proof of condition 4, $\zeta''(s)\neq0 \quad (\sigma<\frac{1}{2},\, t\geq t_5-1)$.
Thus, assuming RH, we only have to consider for $\sigma=\frac{1}{2}$ and $\sigma\in\left[\frac{1}{2},5\right]$ for $\zeta(s)$ and $\zeta''(s)$, respectively.
We take $t_0\in[t_5+1, t_5+2]$ for which
\begin{equation*}
\zeta\left(\frac{1}{2}+it_0\right)\neq0 \quad \text{and} \quad \zeta''(\sigma+it_0)\neq0,\,\, \sigma\in\left[\frac{1}{2},5\right]
\end{equation*}
hold.
Note that this is possible by the identity theorem for complex analytic functions.
Then, we have shown that $t_0$ defined above satisfies $t_0\geq\max\{30,-\sigma_0\}$ and also conditions 3 to 5.\\
\end{proof}

\vspace{3mm}
\begin{prop} \label{prop2}
Assume RH. Take $a$ and $t_0$ which satisfy all conditions of Lemma \ref{lem1}.
Then for $T\geq t_0$ which satisfies $\zeta''(\sigma+iT)\neq0$ and $\zeta(\sigma+iT)\neq0$ for any $\sigma\in\mathbb{R}$, we have
\begin{equation*}
\begin{aligned}
\sum_{\substack{\rho''=\beta''+i\gamma'',\\ 0<\gamma''\leq T}} \left(\beta''-\frac{1}{2}\right) &= \frac{2T}{2\pi}\log{\log{\frac{T}{2\pi}}} + \frac{1}{2\pi}\left(\frac{1}{2}\log{2} - 2\log{\log{2}}\right)T - 2\textnormal{Li}\left(\frac{T}{2\pi}\right)\\
&\quad\quad\quad+ \frac{1}{2\pi}\int_{1/2}^a \left(-\arg{\zeta(\sigma+iT)}+\arg{G_2(\sigma+iT)}\right)d\sigma + O_{a,t_0}(1)
\end{aligned}
\end{equation*}
where the logarithmic branches are taken so that $\log{\zeta(s)}$ and $\log{G_2(s)}$ tend to 0 as $\sigma\rightarrow\infty$ and are holomorphic in $\mathbb{C}\backslash\{\rho+\lambda \,|\, \zeta(\rho)=0 \,\,\text{or}\,\, \infty, \, \lambda\leq0\}$ and $\mathbb{C}\backslash\{\rho''+\lambda \,|\, \zeta''(\rho'')=0 \,\,\text{or}\,\, \infty, \, \lambda\leq0\}$, respectively.
\end{prop}

\vspace{1mm}
\begin{proof}
The steps of the proof generally follow the proof of Proposition 2.2 of \cite{aka}.
We first take $a$, $\sigma_0$, and $t_0$ as in Lemma \ref{lem1} and fix them.
Then, we take $T\geq t_0$ such that $\zeta''(\sigma+iT)\neq0$ and $\zeta(\sigma+iT)\neq0$ ($^\forall\sigma\in\mathbb{R}$).
We also let $\delta\in(0,\frac{1}{2}]$ and put $b:=\frac{1}{2}-\delta$.
We consider the rectangle with vertices $b+it_0$,  $a+it_0$, $a+iT$, and $b+iT$, and then we apply Littlewood's lemma (cf. \cite[pp. 132--133]{tit1}) to $G_2(s)$ there.
By taking the imaginary part, we obtain
\begin{equation} \label{eq:lit}
\begin{aligned}
2\pi\sum_{\substack{\rho''=\beta''+i\gamma'',\\ t_0<\gamma''\leq T}} (\beta''-b) &=
\int_{t_0}^T\log{|G_2(b+it)|}dt - \int_{t_0}^T\log{|G_2(a+it)|}dt\\
&\quad\quad\quad\quad\quad\quad\quad- \int_b^a\arg{G_2(\sigma+it_0)}d\sigma + \int_b^a\arg{G_2(\sigma+iT)}d\sigma\\
&=: I_1 + I_2 + I_3 + \int_b^a\arg{G_2(\sigma+iT)}d\sigma
\end{aligned}
\end{equation}
where the sum is counted with multiplicity.
By the same reasoning as in \cite[p. 2246]{aka}, we have
\begin{equation*}
I_2 = O_a(1),\quad I_3 = O_{a,t_0}(1).
\end{equation*}

Now we only need to estimate $I_1$.
From the functional equation for $\zeta(s)$ (eq. \eqref{eq:fe}), we can deduce that
\begin{equation*}
\begin{aligned}
\zeta''(s) &= F''(s)\zeta(1-s)\left\{1 - 2\frac{1}{\frac{F''}{F'}(s)}\frac{\zeta'}{\zeta}(1-s) + \frac{1}{\frac{F''}{F}(s)}\frac{\zeta''}{\zeta}(1-s)\right\}\\
&= F(s)\frac{F''}{F}(s)\zeta(1-s)\left\{1 - 2\frac{1}{\frac{F''}{F'}(s)}\frac{\zeta'}{\zeta}(1-s) + \frac{1}{\frac{F''}{F}(s)}\frac{\zeta''}{\zeta}(1-s)\right\}.
\end{aligned}
\end{equation*}
Hence,
\begin{equation} \label{eq:i1}
\begin{aligned}
I_1 &= \int_{t_0}^T\log{|G_2(b+it)|}dt = \int_{t_0}^T\log{\frac{2^b}{(\log{2})^2}|\zeta''(b+it)|}dt\\
&= \int_{t_0}^T\log{\frac{2^b}{(\log{2})^2}}dt +  \int_{t_0}^T\log{|\zeta''(b+it)|}dt\\
&= (b\log{2}-2\log{\log{2}})(T-t_0) + \int_{t_0}^T\log{|F(b+it)|}dt + \int_{t_0}^T\log{\left|\frac{F''}{F}(b+it)\right|}dt\quad\quad\quad\\
&\quad\quad\quad\quad\quad+ \int_{t_0}^T\log{|\zeta(1-b-it)|}dt\\
&\quad\quad\quad\quad\quad+ \int_{t_0}^T\log{\left|1 - 2\frac{1}{\frac{F''}{F'}(b+it)}\frac{\zeta'}{\zeta}(1-b-it) + \frac{1}{\frac{F''}{F}(b+it)}\frac{\zeta''}{\zeta}(1-b-it)\right|}dt\\
&=: ((b\log{2}-2\log{\log{2}})T + O_{t_0}(1)) + I_{12} + I_{13} + I_{14} + I_{15}.
\end{aligned}
\end{equation}

Referring to \cite[pp. 2247--2249]{aka}, we have
\begin{equation*}
I_{12} = \left(\frac{1}{2}-b\right)\left(T\log{\frac{T}{2\pi}}-T\right) + O_{t_0}(1),
\end{equation*}
\begin{equation*}
I_{14} = -\int_{1-b}^a\arg{\zeta(\sigma+iT)}d\sigma + O_{a,t_0}(1).
\end{equation*}

Below, we estimate $I_{13}$ and $I_{15}$. We begin with the estimation of $I_{13}$.
As what we did in the proof of condition 2 of Lemma \ref{lem1}, we use Stirling's formula to $\frac{\Gamma'}{\Gamma}(1-s)$ in the region $-\frac{\pi}{2}\leq\arg{(1-s)}\leq\frac{\pi}{2}$ and product rule of derivatives for $\frac{F'}{F}(s)$ (eq. \eqref{eq:f2}), and we can show that
\begin{equation*}
\begin{aligned}
\frac{F''}{F}(b+it) &= \left(\log{\frac{t}{2\pi}}\right)^2 \times
\Bigg\{1-\frac{(1-b)t\log{\frac{t}{2\pi}}-2(1-b)t}{((1-b)^2+t^2)t\log^2{\frac{t}{2\pi}}}\\
&\quad\quad-\frac{i\left\{(1-b)^2-t^2-2(1-b)^3\log{\frac{t}{2\pi}}-(1-2b)t^2\log{\frac{t}{2\pi}}\right\}}{((1-b)^2+t^2)t\log^2{\frac{t}{2\pi}}} + O\left(\frac{1}{t^2\log{t}}\right)\Bigg\}.
\end{aligned}
\end{equation*}
Therefore,
\begin{equation*}
\begin{aligned}
\text{Re}\left(\log{\frac{F''}{F}(b+it)}\right) &= 2\log{\log{\frac{t}{2\pi}}} -\frac{(1-b)t\log{\frac{t}{2\pi}}}{((1-b)^2+t^2)t\log^2{\frac{t}{2\pi}}} + O\left(\frac{1}{t^2\log{t}}\right)\\
&= 2\log{\log{\frac{t}{2\pi}}} + O\left(\frac{1}{t^2\log{t}}\right).
\end{aligned}
\end{equation*}
Hence,
\begin{equation*}
\begin{aligned}
I_{13} &= \int_{t_0}^T\log{\left|\frac{F''}{F}(b+it)\right|}dt =  \int_{t_0}^T\text{Re}\left(\log{\frac{F''}{F}(b+it)}\right)dt\\
&= 2\int_{t_0}^T \log{\log{\frac{t}{2\pi}}}dt + O\left(\int_{t_0}^T\frac{dt}{t^2\log{t}}\right)\\
&= 2T\log{\log{\frac{T}{2\pi}}} - 4\pi\text{Li}\left(\frac{T}{2\pi}\right) + O_{t_0}(1).
\end{aligned}
\end{equation*}

Finally, we estimate $I_{15}$.
Again from the functional equation for $\zeta(s)$ (eq. \eqref{eq:fe}), we have
\begin{equation} \label{eq:fe2}
\frac{1}{\frac{F''}{F}(s)}\frac{\zeta''}{\zeta}(s) = 1 - \left\{2\frac{1}{\frac{F''}{F'}(s)}\frac{\zeta'}{\zeta}(1-s) - \frac{1}{\frac{F''}{F}(s)}\frac{\zeta''}{\zeta}(1-s)\right\}.
\end{equation}
It follows from condition 2 of Lemma \ref{lem1} that the right hand side of equation \eqref{eq:fe2} is holomorphic and has no zeros in the region defined by $\sigma\leq\sigma_0$ and $t\geq2$.
And from conditions 3 and 4 of Lemma \ref{lem1}, the left hand side of equation \eqref{eq:fe2} is holomorphic and has no zeros in the region defined by $\sigma_0\leq\sigma<\frac{1}{2}$ and $t\geq t_0-1$.
Thus, we can determine $\log{\left(1 - \left\{2\frac{1}{\frac{F''}{F'}(s)}\frac{\zeta'}{\zeta}(1-s) - \frac{1}{\frac{F''}{F}(s)}\frac{\zeta''}{\zeta}(1-s)\right\}\right)}$ so that it tends to $0$ as $\sigma\rightarrow-\infty$ which follows from condition 2 of Lemma \ref{lem1}, and is holomorphic in the region $\sigma<\frac{1}{2},\, t>t_0-1$.
Now we consider the same trapezoid $C$ as in \cite[p. 2247]{aka}.
Then by Cauchy's integral theorem,
\begin{equation} \label{eq:i15}
\int_C \log{\left(1 - 2\frac{1}{\frac{F''}{F'}(s)}\frac{\zeta'}{\zeta}(1-s) + \frac{1}{\frac{F''}{F}(s)}\frac{\zeta''}{\zeta}(1-s)\right)}ds =0.
\end{equation}
By using condition 2 of Lemma \ref{lem1}, we can also show that (cf. \cite[p. 2248]{aka})
\begin{equation*}
\left|\left(\int_{\sigma_0+iT}^{-T+iT}+\int_{-T+iT}^{-t_0+it_0}+\int_{-t_0+it_0}^{\sigma_0+it_0}\right)\log{\left(1 - 2\frac{1}{\frac{F''}{F'}(s)}\frac{\zeta'}{\zeta}(1-s) + \frac{1}{\frac{F''}{F}(s)}\frac{\zeta''}{\zeta}(1-s)\right)}ds\right| \ll 1.
\end{equation*}

Next, we estimate the integral from $\sigma_0+it_0$ to $b+it_0$ trivially, we obtain
\begin{equation*}
\int_{\sigma_0+it_0}^{b+it_0} \log{\left(1 - 2\frac{1}{\frac{F''}{F'}(s)}\frac{\zeta'}{\zeta}(1-s) + \frac{1}{\frac{F''}{F}(s)}\frac{\zeta''}{\zeta}(1-s)\right)} ds = O_{t_0}(1).
\end{equation*}

Substituting the above two equations into equation \eqref{eq:i15} and taking the imaginary part, we obtain
\begin{equation*}
\begin{aligned}
&I_{15} =  \int_{t_0}^T \log{\left|1 - 2\frac{1}{\frac{F''}{F'}(b+it)}\frac{\zeta'}{\zeta}(1-b-it) + \frac{1}{\frac{F''}{F}(b+it)}\frac{\zeta''}{\zeta}(1-b-it)\right|} dt\\
&= \int_{\sigma_0}^b \arg{\left(1 - 2\frac{1}{\frac{F''}{F'}(\sigma+iT)}\frac{\zeta'}{\zeta}(1-\sigma-iT) + \frac{1}{\frac{F''}{F}(\sigma+iT)}\frac{\zeta''}{\zeta}(1-\sigma-iT)\right)} d\sigma + O_{t_0}(1)\\
&\overset{\eqref{eq:fe2}}= \int_{\sigma_0}^b \arg{\left(\frac{1}{\frac{F''}{F}(\sigma+iT)}\frac{\zeta''}{\zeta}(\sigma+iT)\right)}d\sigma + O_{t_0}(1).
\end{aligned}
\end{equation*}
From conditions 3 and 4 of Lemma \ref{lem1}, there exists an $n\in\mathbb{Z}$ such that
\begin{equation*}
\left(2n-\frac{7}{6}\right)\pi \leq \arg{\left(\frac{1}{\frac{F''}{F}(\sigma+iT)}\frac{\zeta''}{\zeta}(\sigma+iT)\right)} \leq \left(2n+\frac{7}{6}\right)\pi
\end{equation*}
for $\sigma_0\leq\sigma<\frac{1}{2}$.
Again using equation \eqref{eq:fe2}, we note that from condition 2 of Lemma \ref{lem1},
$\text{Re}\left(\frac{1}{\frac{F''}{F}(\sigma_0+iT)}\frac{\zeta''}{\zeta}(\sigma_0+iT)\right)>0$.
From our choice of logarithmic branch, we have
\begin{equation*}
-\frac{\pi}{2}\leq\arg{\left(\frac{1}{\frac{F''}{F}(\sigma_0+iT)}\frac{\zeta''}{\zeta}(\sigma_0+iT)\right)}\leq\frac{\pi}{2}.
\end{equation*}
Now, the only integer $n\in\mathbb{Z}$ such that $[\left(2n-\frac{7}{6}\right)\pi, \left(2n+\frac{7}{6}\right)\pi]$ includes $[-\frac{\pi}{2}, \frac{\pi}{2}]$ is $n=0$.
Hence for $\sigma_0\leq\sigma<\frac{1}{2}$, we have
\begin{equation} \label{eq:arg}
-\frac{7\pi}{6} \leq \arg{\left(\frac{1}{\frac{F''}{F}(\sigma+iT)}\frac{\zeta''}{\zeta}(\sigma+iT)\right)} \leq \frac{7\pi}{6}
\end{equation}
which gives us immediately
\begin{equation*}
I_{15} = O_{t_0}(1).
\end{equation*}

Substituting the estimations on $I_{12}$, $I_{13}$, $I_{14}$, and $I_{15}$ into equation \eqref{eq:i1}, we obtain
\begin{equation*}
\begin{aligned}
I_1 = (b\log{2}-2\log{\log{2}})T &+ \left(\frac{1}{2}-b\right)\left(T\log{\frac{T}{2\pi}}-T\right) + 2T\log{\log{\frac{T}{2\pi}}} - 4\pi\text{Li}\left(\frac{T}{2\pi}\right)\\
&\quad\quad\quad\quad\quad\quad\quad\quad\quad\quad\quad- \int_{1-b}^a\arg{\zeta(\sigma+iT)}d\sigma + O_{a,t_0}(1).
\end{aligned}
\end{equation*}

To finalize the proof of Proposition \ref{prop2}, we substitute the estimations on $I_1$, $I_2$, and $I_3$ into equation \eqref{eq:lit} to obtain
\begin{equation*}
\begin{aligned}
2&\pi\sum_{\substack{\rho''=\beta''+i\gamma'',\\ 0<\gamma''\leq T}} (\beta''-b) =
2T\log{\log{\frac{T}{2\pi}}} + (b\log{2}-2\log{\log{2}})T - 4\pi\text{Li}\left(\frac{T}{2\pi}\right)\\
&+\left(\frac{1}{2}-b\right)\left(T\log{\frac{T}{2\pi}}-T\right) - \int_{1-b}^a\arg{\zeta(\sigma+iT)}d\sigma + \int_b^a\arg{G_2(\sigma+iT)}d\sigma + O_{a,t_0}(1).
\end{aligned}
\end{equation*}
Taking the limit $\delta\rightarrow0$, we have $b\rightarrow\frac{1}{2}$, thus
\begin{equation*}
\begin{aligned}
\sum_{\substack{\rho''=\beta''+i\gamma'',\\ 0<\gamma''\leq T}} \left(\beta''-\frac{1}{2}\right) &= \frac{2T}{2\pi}\log{\log{\frac{T}{2\pi}}} + \frac{1}{2\pi}\left(\frac{1}{2}\log{2}-2\log{\log{2}}\right)T - 2\text{Li}\left(\frac{T}{2\pi}\right)\\
&\quad\quad\quad\quad+ \frac{1}{2\pi}\int_{\frac{1}{2}}^a(-\arg{\zeta(\sigma+iT)}+\arg{G_2(\sigma+iT)})d\sigma + O_{a,t_0}(1).
\end{aligned}
\end{equation*}
\end{proof}

\vspace{3mm}
To complete the proof of Theorem \ref{cha1}, we need to estimate
\begin{equation*}
\int_{\frac{1}{2}}^a(-\arg{\zeta(\sigma+iT)}+\arg{G_2(\sigma+iT)})d\sigma
\end{equation*}
in Proposition \ref{prop2}.
For that purpose, similar to the method taken in \cite{aka}, below we give two bounds for $-\arg{\zeta(\sigma+iT)}+\arg{G_2(\sigma+iT)}$.
We write
\begin{equation*}
-\arg{\zeta(\sigma+iT)}+\arg{G_2(\sigma+iT)} = \arg{\frac{G_2}{\zeta}(\sigma+iT)}
\end{equation*}
and take the argument on the right hand side so that
$\log{\frac{G_2}{\zeta}(\sigma+iT)}$ tends to $0$ as $\sigma\rightarrow\infty$ and is holomorphic in $\mathbb{C}\backslash\{z+\lambda \,|\, \frac{\zeta''}{\zeta}(z)=0 \,\,\text{or}\,\, \infty, \, \lambda\leq0\}$.

\vspace{3mm}
\begin{lem} \label{lem3}
Assume RH and let $T\geq t_0$.
Then for any $\epsilon_0>0$ satisfying $\epsilon_0<\frac{3}{8}\frac{1}{\log{T}}$
(since $T\geq t_0\geq30$, $\epsilon_0<\frac{1}{8}$),
we have for $\frac{1}{2}+\epsilon_0<\sigma\leq a$,
\begin{equation*}
\arg{\frac{G_2}{\zeta}(\sigma+iT)} = O_{a,t_0} \left(\frac{\log{\frac{\log{T}}{\epsilon_0}}}{\sigma-\frac{1}{2}-\epsilon_0}\right).
\end{equation*}
\end{lem}

\vspace{1mm}
\begin{proof}
To begin, we note that $\frac{G_2}{\zeta}(s)$ is uniformly convergent to $1$ as $\sigma\rightarrow\infty$ for $t\in\mathbb{R}$, so we can take a number $c\in\mathbb{R}$ satisfying $a+1\leq c\leq\frac{t_0}{2}$ and $\frac{1}{2}\leq\text{Re}\left(\frac{G_2}{\zeta}(s)\right)\leq\frac{3}{2}\quad (\sigma\geq c)$. In fact, we can check that taking $c=15$ is enough.
Now, the proof also proceeds similarly to the proof of Lemma 2.3 of \cite{aka}.
We let $\sigma\in(\frac{1}{2}+\epsilon_0,a]$ and let
$q_{G_2/\zeta}=q_{G_2/\zeta}(\sigma,T)$ denote the number of times $\text{Re}\left(\frac{G_2}{\zeta}(u+iT)\right)$ vanishes in $u\in[\sigma,c]$.
Then, $\left|\arg{\frac{G_2}{\zeta}(\sigma+iT)}\right|\leq\left(q_{G_2/\zeta}+1\right)\pi$.
Now we estimate $q_{G_2/\zeta}$.
For that purpose, we set
\begin{equation*}
H_2(z)=H_{2_T}(z):=\frac{\frac{G_2}{\zeta}(z+iT)+\frac{G_2}{\zeta}(z-iT)}{2} \quad (z\in\mathbb{C})
\end{equation*}
and $n_{H_2}(r):=\sharp\{z\in\mathbb{C} \,|\, H_2(z)=0, |z-c|\leq r\}$.
Then, we have $q_{G_2/\zeta}\leq n_{H_2}(c-\sigma)$ for $\frac{1}{2}+\epsilon_0<\sigma\leq a$.
For each $\sigma\in(\frac{1}{2}+\epsilon_0,a]$, we take $\epsilon=\epsilon_{\sigma,T}$ satisfying $0<\epsilon<\sigma-\frac{1}{2}-\epsilon_0$, then $H_2(z)$ is holomorphic in the region $\{z\in\mathbb{C} \,|\, |z-c| \leq c-\sigma+\epsilon\}$.
As in \cite[p. 2250]{aka}, by using Jensen's theorem (cf. \cite[pp. 125--126]{tit1}), we can show that
\begin{equation*}
\begin{aligned}
n_{H_2}(c-\sigma) &\leq \frac{1}{C_1\epsilon}\int_0^{c-\sigma+\epsilon}\frac{n_{H_2}(r)}{r}dr\\
&= \frac{1}{C_1\epsilon}\frac{1}{2\pi}\int_0^{2\pi}\log{|H_2(c+(c-\sigma+\epsilon)e^{i\theta})|}d\theta - \frac{1}{C_1\epsilon}\log{|H_2(c)|}
\end{aligned}
\end{equation*}
for some constant $C_1>0$,
which by our choice of $c$ gives us
\begin{equation} \label{eq:nh2}
n_{H_2}(c-\sigma) \leq \frac{1}{C_1\epsilon}\frac{1}{2\pi}\int_0^{2\pi}\log{|H_2(c+(c-\sigma+\epsilon)e^{i\theta})|}d\theta + \frac{1}{\epsilon}O_{a,t_0}(1).
\end{equation}

As noted in \cite[p. 2250]{aka} (cf. \cite[Theorems 9.2 and 9.6(A)]{tit2}),
\begin{equation*}
\frac{\zeta'}{\zeta}(\sigma\pm it) = O\left(\frac{\log{T}}{\sigma-\frac{1}{2}}\right)
\end{equation*}
holds for $\frac{T}{2}\leq t\leq 2T$ and $\frac{1}{2}<\sigma\leq 2c$.
Thus, for $\frac{T}{2}\leq t\leq 2T$ and $\frac{1}{2}+\epsilon_0<\sigma\leq 2c$, we have
$\frac{\zeta'}{\zeta}(\sigma\pm it) = O\left(\frac{\log{T}}{\epsilon_0}\right)$, and
\begin{equation*}
\left(\frac{\zeta'}{\zeta}(\sigma\pm it)\right)^2 = O\left(\frac{\log^2{T}}{\epsilon_0^2}\right), \quad\quad \text{for}\quad \frac{T}{2}\leq t\leq 2T,\, \frac{1}{2}+\epsilon_0<\sigma\leq 2c
\end{equation*}
immediately follows.
Applying Cauchy's integral formula, we have for $s=\sigma+it$ with $\frac{T}{2}<|t|< 2T$ and $\frac{1}{2}+\epsilon_0<\sigma<2c$,
\begin{equation*}
\left(\frac{\zeta'}{\zeta}(s)\right)' = \frac{1}{2\pi i}\int_{|z-s|=\epsilon_0} \frac{\frac{\zeta'}{\zeta}(z)}{(z-s)^2}dz = O\left(\frac{\log{T}}{\epsilon_0^2}\right).
\end{equation*}

Therefore,
\begin{equation*}
\frac{\zeta''}{\zeta}(s) = \left(\frac{\zeta'}{\zeta}(s)\right)' + \left(\frac{\zeta'}{\zeta}(s)\right)^2 = O\left(\frac{\log^2{T}}{\epsilon_0^2}\right)
\end{equation*}
for $\frac{T}{2}<|t|< 2T$ and $\frac{1}{2}+\epsilon_0<\sigma<2c$.\\

Consequently, we can easily show that
\begin{equation*}
|H_2(c+(c-\sigma+\epsilon)e^{i\theta})| \ll_a \frac{\log^2{T}}{\epsilon_0^2}
\end{equation*}
and so
\begin{equation*}
\frac{1}{2\pi}\int_0^{2\pi}\log{|H_2(c+(c-\sigma+\epsilon)e^{i\theta})|}d\theta \ll_a \log{\frac{\log{T}}{\epsilon_0}}.
\end{equation*}
Substituting this into equation \eqref{eq:nh2}, we obtain
\begin{equation*}
n_{H_2}(c-\sigma) = \frac{1}{\epsilon}O_{a,t_0}\left(\log{\frac{\log{T}}{\epsilon_0}}\right)
\end{equation*}
which implies
\begin{equation*}
\arg{\frac{G_2}{\zeta}(\sigma+iT)} = \frac{1}{\epsilon}O_{a,t_0}\left(\log{\frac{\log{T}}{\epsilon_0}}\right).
\end{equation*}
Taking $\epsilon=\frac{1}{2}\left(\sigma-\frac{1}{2}-\epsilon_0\right) \,\, \left(<\sigma-\frac{1}{2}-\epsilon_0\right)$, we obtain
\begin{equation*}
\arg{\frac{G_2}{\zeta}(\sigma+iT)} = O_{a,t_0}\left(\frac{\log{\frac{\log{T}}{\epsilon_0}}}{\sigma-\frac{1}{2}-\epsilon_0}\right).
\end{equation*}
\end{proof}

\vspace{3mm}
\begin{lem} \label{lem4}
Assume RH and let $A>3$ be fixed. Then there exists a constant $C_0>0$ such that
\begin{equation*}
\left|\zeta''(\sigma+it)\right|\leq\exp{\left( C_0\left(\frac{(\log T)^{2(1-\sigma)}}{\log\log{T}}+(\log{T})^{1/10}\right)\right)}
\end{equation*}
holds for $T\geq30$, $\frac{T}{2}\leq t\leq 2T$, $\frac{1}{2}-\frac{1}{\log{\log{T}}}\leq\sigma\leq A$.
\end{lem}

\vspace{1mm}
\begin{proof}
Referring to \cite[equations (14.14.2), (14.14.5) and the first equation on p. 384]{tit2} (cf. \cite[pp. 2251--2252]{aka}), we know that
\begin{equation*}
|\zeta(\sigma+it)| \leq \exp{\left(C_2\left(\frac{(\log{T})^{2(1-\sigma)}}{\log{\log{T}}}\right)+(\log{T})^{1/10}\right)}
\end{equation*}
holds for  $\frac{1}{2}-\frac{2}{\log{\log{T}}}\leq\sigma\leq A+1$, $\frac{T}{3}\leq t\leq 3T$ for some constant $C_2>0$.\\

Applying Cauchy's integral formula,
\begin{equation*}
\zeta''(s) = \frac{1}{\pi i} \int_{|z-s|=\epsilon} \frac{\zeta(z)}{(z-s)^3}dz \quad\quad \text{for} \quad 0<\epsilon<\frac{1}{2}
\end{equation*}
in the region defined by  $\frac{1}{2}-\frac{1}{\log{\log{T}}}\leq\sigma\leq A$ and $\frac{T}{2}\leq t\leq 2T$.
Thus, in that region we have
\begin{equation*}
|\zeta''(s)| \leq \frac{2}{\epsilon^2} \exp{\left(C_2\left(\frac{(\log{T})^{2(1-\sigma)}}{\log{\log{T}}}\right)+(\log{T})^{1/10}\right)}.
\end{equation*}
Taking $\epsilon=\frac{1}{2\sqrt{\log{\log{T}}}} \,\, (<\frac{1}{2})$, we obtain Lemma \ref{lem4}.\\
\end{proof}

\vspace{3mm}
\begin{lem} \label{lem5}
Assume RH and let $T\geq30$.
Then for any $\frac{1}{2}\leq\sigma\leq\frac{3}{4}$, we have
\begin{equation*}
\arg{G_2(\sigma+iT)} = O_a\left(\frac{(\log T)^{2(1-\sigma)}}{(\log{\log{T}})^{1/2}}\right).
\end{equation*}
\end{lem}

\vspace{1mm}
\begin{proof}
The proof proceeds in the same way as the proof of Lemma 2.4 of \cite{aka}.
Refer to \cite[pp. 2252--2253]{aka} for the detailed proof and use Lemma \ref{lem4} above in place of Lemma 2.6 of \cite{aka}.\\
\end{proof}

\vspace{3mm}
\begin{rmk}
The restrictions of the lower bound of $T$ we gave in Lemmas \ref{lem3}, \ref{lem4}, and \ref{lem5} are in fact not essential.
Nevertheless, they are sufficient for our needs.

We may let $T$ to be any positive number in Lemmas \ref{lem3}, \ref{lem4}, and \ref{lem5}, however in that case, we need to modify some calculations in the proofs. Thus, we used these restrictions for our convenience.
\end{rmk}

\vspace{3mm}
\noindent
\underline{\textbf{Proof of Theorem \ref{cha1}:}}\\

First of all, we consider for $T\geq t_0$ which satisfies $\zeta''(\sigma+iT)\neq0$ and $\zeta(\sigma+iT)\neq0$ for any $\sigma\in\mathbb{R}$. By Lemma \ref{lem3}, we have
\begin{equation*}
\int_{\frac{1}{2}+2\epsilon_0}^a \arg{\frac{G_2}{\zeta}(\sigma+iT)}d\sigma \ll_{a,t_0} \int_{\frac{1}{2}+2\epsilon_0}^a \frac{\log{\frac{\log{T}}{\epsilon_0}}}{\sigma-\frac{1}{2}-\epsilon_0}d\sigma
\ll_a \log{\frac{\log{T}}{\epsilon_0}}\log{\frac{1}{\epsilon_0}}.
\end{equation*}
Next, by Lemma \ref{lem5},
\begin{equation*}
\arg{G_2(\sigma+iT)} = O_a\left(\frac{(\log T)^{2(1-\sigma)}}{(\log{\log{T}})^{1/2}}\right) \quad\quad \text{for} \quad \frac{1}{2}\leq\sigma\leq\frac{3}{4}
\end{equation*}
and from equation (2.23) of \cite[p. 2251]{aka} (cf. \cite[equations (14.14.3) and (14.14.5)]{tit2}),
RH implies that
\begin{equation*}
\arg{\zeta(\sigma+iT)} = O\left(\frac{(\log T)^{2(1-\sigma)}}{\log{\log{T}}}\right)
\end{equation*}
holds uniformly for $\frac{1}{2}\leq\sigma\leq\frac{3}{4}$.
Thus,
\begin{equation*}
\int_{\frac{1}{2}}^{\frac{1}{2}+2\epsilon_0} \arg{\frac{G_2}{\zeta}(\sigma+iT)}d\sigma \ll_a \frac{\log T}{(\log{\log{T}})^{1/2}}\epsilon_0.
\end{equation*}
Now we take $\epsilon_0 = \frac{1}{4\log{T}} \,\, \left(<\frac{3}{8}\frac{1}{\log{T}}\right)$, then we have
\begin{equation*}
\int_{\frac{1}{2}}^a \arg{\frac{G_2}{\zeta}(\sigma+iT)}d\sigma \ll_{a,t_0} (\log{\log{T}})^2.
\end{equation*}
Applying this to Proposition \ref{prop2}, we have
\begin{equation} \label{eq:kari1}
\begin{aligned}
\sum_{\substack{\rho''=\beta''+i\gamma'',\\ 0<\gamma''\leq T}} \left(\beta''-\frac{1}{2}\right) = \,\frac{2T}{2\pi}\log{\log{\frac{T}{2\pi}}} + \frac{1}{2\pi}\left(\frac{1}{2}\log{2} - 2\log{\log{2}}\right)&T - 2\textnormal{Li}\left(\frac{T}{2\pi}\right)\\
&+ O_{a,t_0}((\log{\log{T}})^2).
\end{aligned}
\end{equation}

Secondly, for $2\pi<T<t_0$, we are adding some finite number of terms which depend on $t_0$ so this can be included in the error term.\\

Thirdly, for $T\geq t_0$ such that $\zeta''(\sigma+iT)=0$ or $\zeta(\sigma+iT)=0$ for some $\sigma\in\mathbb{R}$, then there is some increment in the value of $\sum_{\substack{\rho''=\beta''+i\gamma'',\\ 0<\gamma''\leq T}} \left(\beta''-\frac{1}{2}\right)$ as much as $\sum_{\substack{\rho''=\beta''+i\gamma'',\\ \gamma'' = T}} \left(\beta''-\frac{1}{2}\right)$. Now we estimate this and we show that this can be included in the error term of equation \eqref{eq:kari1}.
We start by taking a small $0<\epsilon<1$ such that $\zeta''(\sigma+i(T\pm\epsilon))\neq0$ and $\zeta(\sigma+i(T\pm\epsilon))\neq0$ for any $\sigma\in\mathbb{R}$.
According to equation \eqref{eq:kari1},
\begin{equation*}
\begin{aligned}
\sum_{\substack{\rho''=\beta''+i\gamma'',\\ 0<\gamma''\leq T+\epsilon}} \left(\beta''-\frac{1}{2}\right) = \frac{2(T+\epsilon)}{2\pi}\log{\log{\frac{T+\epsilon}{2\pi}}} &+ \frac{1}{2\pi}\left(\frac{1}{2}\log{2} - 2\log{\log{2}}\right)(T+\epsilon)\\
&- 2\textnormal{Li}\left(\frac{T+\epsilon}{2\pi}\right) + O_{a,t_0}((\log{\log{T}})^2),
\end{aligned}
\end{equation*}
\begin{equation*}
\begin{aligned}
\sum_{\substack{\rho''=\beta''+i\gamma'',\\ 0<\gamma''\leq T-\epsilon}} \left(\beta''-\frac{1}{2}\right) = \frac{2(T-\epsilon)}{2\pi}\log{\log{\frac{T-\epsilon}{2\pi}}} &+ \frac{1}{2\pi}\left(\frac{1}{2}\log{2} - 2\log{\log{2}}\right)(T-\epsilon)\\
&- 2\textnormal{Li}\left(\frac{T-\epsilon}{2\pi}\right) + O_{a,t_0}((\log{\log{T}})^2).
\end{aligned}
\end{equation*}
Thus,
\begin{equation*}
\begin{aligned}
\sum_{\substack{\rho''=\beta''+i\gamma'',\\ T-\epsilon<\gamma''\leq T+\epsilon}} \left(\beta''-\frac{1}{2}\right) &= \frac{2(T+\epsilon)}{2\pi}\log{\log{\frac{T+\epsilon}{2\pi}}} - \frac{2(T-\epsilon)}{2\pi}\log{\log{\frac{T-\epsilon}{2\pi}}}\\
&\quad\quad+ \frac{\epsilon}{\pi}\left(\frac{1}{2}\log{2} - 2\log{\log{2}}\right) - 2\left\{\textnormal{Li}\left(\frac{T+\epsilon}{2\pi}\right) - \textnormal{Li}\left(\frac{T-\epsilon}{2\pi}\right)\right\}\\
&\quad\quad\quad\quad\quad\quad\quad\quad\quad\quad\quad\quad\quad\quad\quad\quad\quad\quad\quad\quad+ O_{a,t_0}((\log{\log{T}})^2)\\
&= \frac{2\epsilon}{\pi}\log{\log{\frac{T}{2\pi}}} + \frac{2\epsilon}{\pi\log{\frac{T}{2\pi}}} + \frac{\epsilon}{\pi}\left(\frac{1}{2}\log{2} - 2\log{\log{2}}\right) - \frac{2\epsilon}{\pi\log{\frac{T}{2\pi}}}\\
&\,\,\quad\quad\quad\quad\quad\quad\quad\quad\quad\quad\quad\quad+ O\left(\frac{\epsilon^2}{T\log{T}}\right) + O_{a,t_0}((\log{\log{T}})^2).
\end{aligned}
\end{equation*}
This gives us
\begin{equation*}
\sum_{\substack{\rho''=\beta''+i\gamma'',\\ T-\epsilon<\gamma''\leq T+\epsilon}} \left(\beta''-\frac{1}{2}\right) = O_{a,t_0}((\log{\log{T}})^2)
\end{equation*}
which implies
\begin{equation*}
\sum_{\substack{\rho''=\beta''+i\gamma'',\\ \gamma'' = T}} \left(\beta''-\frac{1}{2}\right) = O_{a,t_0}((\log{\log{T}})^2).
\end{equation*}
Therefore, this increment can also be included in the error term.\\

Finally, since $a$ and $t_0$ are fixed constants,
\begin{equation*}
\begin{aligned}
\sum_{\substack{\rho''=\beta''+i\gamma'',\\ 0<\gamma''\leq T}} \left(\beta''-\frac{1}{2}\right) = \,\frac{2T}{2\pi}\log{\log{\frac{T}{2\pi}}} + \frac{1}{2\pi}\left(\frac{1}{2}\log{2} - 2\log{\log{2}}\right)T &- 2\textnormal{Li}\left(\frac{T}{2\pi}\right)\\
&+ O((\log{\log{T}})^2)
\end{aligned}
\end{equation*}
holds for any $T>2\pi$.\\
\qed

\vspace{3mm}
\noindent
\underline{\textbf{Proof of Corollary \ref{cha12}:}}\\

This is an immediate consequence of Theorem \ref{cha1}.
For the proof, refer to \cite[p. 58 (ending part of section 3)]{lev}.\\
\qed


\section{Proof of Theorem \ref{cha2}}
\label{sec:3}

In this section we give the proof of Theorem \ref{cha2}.
We first show the following proposition.

\vspace{3mm}
\begin{prop} \label{prop6}
Assume RH. Then for $T\geq2$ which satisfies $\zeta(\sigma+iT)\neq0$ and $\zeta''(\sigma+iT)\neq0$ $(^\forall\sigma\in\mathbb{R})$, we have
\begin{equation*}
N_2(T) = \frac{T}{2\pi}\log{\frac{T}{4\pi}}-\frac{T}{2\pi}+\frac{1}{2\pi}\arg{G_2\left(\frac{1}{2}+iT\right)}+\frac{1}{2\pi}\arg{\zeta\left(\frac{1}{2}+iT\right)}+O_{t_0}(1)
\end{equation*}
where the arguments are taken as in Proposition \ref{prop2}.
\end{prop}

\vspace{1mm}
\begin{proof}
The steps of the proof also follow the proof of Proposition 3.1 of \cite{aka}.
We take $a$, $\sigma_0$, $t_0$, $T$, $\delta$, and $b$ as in the beginning of the proof of Proposition \ref{prop2}. We let $b' := \frac{1}{2}-\frac{\delta}{2}$, then replacing $b$ by $b'$ in equation \eqref{eq:lit}, we have

\begin{equation*}
\begin{aligned}
2\pi\sum_{\substack{\rho''=\beta''+i\gamma'',\\ 0<\gamma''\leq T}} (\beta''-b') &=
\int_{t_0}^T\log{|G_2(b'+it)|}dt - \int_{t_0}^T\log{|G_2(a+it)|}dt\\
&\quad\quad\quad\quad\quad\quad\quad- \int_{b'}^a\arg{G_2(\sigma+it_0)}d\sigma + \int_{b'}^a\arg{G_2(\sigma+iT)}d\sigma.
\end{aligned}
\end{equation*}
Subtracting this from equation \eqref{eq:lit}, we have
\begin{equation} \label{eq:lit'}
\begin{aligned}
\pi\delta (N_2(T) - N_2(t_0)) &=
\int_{t_0}^T\log{|G_2(b+it)|}dt - \int_{t_0}^T\log{|G_2(b'+it)|}dt\\
&\quad\quad\quad\quad\quad\quad- \int_b^{b'}\arg{G_2(\sigma+it_0)}d\sigma + \int_b^{b'}\arg{G_2(\sigma+iT)}d\sigma\\
&=: J_1 + J_2 + J_3 + \int_b^{b'}\arg{G_2(\sigma+iT)}d\sigma.
\end{aligned}
\end{equation}

Referring to the estimation of $I_3$ in the proof of Proposition \ref{prop2} (cf. \cite[p. 2246]{aka}), we can easily show that
\begin{equation*}
J_3 = O_{t_0}(\delta).
\end{equation*}

Now we estimate $J_1 + J_2$.
From equation \eqref{eq:i1}, we have
\begin{equation*}
\begin{aligned}
J_1+J_2 &= \int_{t_0}^T\log{|G_2(b+it)|}dt - \int_{t_0}^T\log{|G_2(b'+it)|}dt\\
&= ((b-b')\log{2})(T-t_0) + \int_{t_0}^T (\log{|F(b+it)|}-\log{|F(b'+it)|})dt\\
&\quad+ \int_{t_0}^T \left(\log{\left|\frac{F''}{F}(b+it)\right|}-\log{\left|\frac{F''}{F}(b'+it)\right|}\right)dt\\
&\quad+ \int_{t_0}^T (\log{|\zeta(1-b-it)|}-\log{|\zeta(1-b'-it)|})dt\\
&\quad+ \int_{t_0}^T \Bigg(\log{\left|1 - 2\frac{1}{\frac{F''}{F'}(b+it)}\frac{\zeta'}{\zeta}(1-b-it) + \frac{1}{\frac{F''}{F}(b+it)}\frac{\zeta''}{\zeta}(1-b-it)\right|}\\
&\quad\quad\quad- \log{\left|1 - 2\frac{1}{\frac{F''}{F'}(b'+it)}\frac{\zeta'}{\zeta}(1-b'-it) + \frac{1}{\frac{F''}{F}(b'+it)}\frac{\zeta''}{\zeta}(1-b'-it)\right|}\Bigg)dt\\
&=: \left(\left(-\frac{\delta}{2}\log{2}\right)T + O_{t_0}(\delta)\right) + J_{12} + J_{13} + J_{14} + J_{15}.
\end{aligned}
\end{equation*}

Referring to \cite[pp. 2255--2256]{aka}, we have
\begin{equation*}
J_{12} = \frac{\delta}{2}\left(T\log{\frac{T}{2\pi}}-T\right) + O_{t_0}(\delta),
\end{equation*}
\begin{equation*}
J_{14} = \int_{1-b'}^{1-b}\arg{\zeta(\sigma+iT)}d\sigma + O_{t_0}(\delta).
\end{equation*}

We only need to estimate $J_{13}$ and $J_{15}$.
We begin with the estimation of $J_{13}$.
Since all zeros and poles of $F(s)$ lie on $\mathbb{R}$, $\frac{F''}{F}(s)$ has no poles in $t>0$.
From condition 3 of Lemma \ref{lem1}, we can define a branch of $\log{\frac{F''}{F}(s)}$ for $0<\sigma<\frac{1}{2}$ and $t>t_0-1$ by $\arg{\frac{F''}{F}(s)}\in\left[-\frac{\pi}{6},\frac{\pi}{6}\right]$.
As in \cite[p. 2255]{aka}, we apply Cauchy's integral theorem to $\log{\frac{F''}{F}(s)}$ on the rectangle with vertices $b+it_0$, $b'+it_0$, $b'+iT$, and $b+iT$ and take the imaginary part, then we obtain
\begin{equation*}
J_{13} =  \int_b^{b'} \arg{\frac{F''}{F}(\sigma+it_0)}d\sigma - \int_b^{b'} \arg{\frac{F''}{F}(\sigma+iT)}d\sigma = O_{t_0}(\delta).
\end{equation*}

Finally, we estimate $J_{15}$. 
We determine a branch of
\begin{equation*}
\log{\left(1 - 2\frac{1}{\frac{F''}{F'}(s)}\frac{\zeta'}{\zeta}(1-s) + \frac{1}{\frac{F''}{F}(s)}\frac{\zeta''}{\zeta}(1-s)\right)}
\end{equation*}
as that in the estimation of $I_{15}$ in the proof of Proposition \ref{prop2}, then it is holomorphic in the region defined by $0<\sigma<\frac{1}{2},\, t>t_0-1$.
Applying Cauchy's integral theorem to it on the path taken for estimating $J_{13}$, we have
\begin{equation*}
\begin{aligned}
J_{15} = &\int_b^{b'} \arg{\left(1 - 2\frac{1}{\frac{F''}{F'}(\sigma+it_0)}\frac{\zeta'}{\zeta}(1-\sigma-it_0) + \frac{1}{\frac{F''}{F}(\sigma+it_0)}\frac{\zeta''}{\zeta}(1-\sigma-it_0)\right)}d\sigma\\
&\quad- \int_b^{b'} \arg{\left(1 - 2\frac{1}{\frac{F''}{F'}(\sigma+iT)}\frac{\zeta'}{\zeta}(1-\sigma-iT) + \frac{1}{\frac{F''}{F}(\sigma+iT)}\frac{\zeta''}{\zeta}(1-\sigma-iT)\right)}d\sigma.
\end{aligned}
\end{equation*}
Again using equation \eqref{eq:fe2}, 
\begin{equation*}
J_{15} = \int_b^{b'} \arg{\left(\frac{1}{\frac{F''}{F}(\sigma+it_0)}\frac{\zeta''}{\zeta}(\sigma+it_0)\right)}d\sigma - \int_b^{b'} \arg{\left(\frac{1}{\frac{F''}{F}(\sigma+iT)}\frac{\zeta''}{\zeta}(\sigma+iT)\right)}d\sigma.
\end{equation*}
Applying inequalities \eqref{eq:arg}, we obtain
\begin{equation*}
J_{15} = O_{t_0}(\delta).
\end{equation*}

In consequence,
\begin{equation*}
J_1+J_2 = \frac{\delta}{2}\left(T\log{\frac{T}{4\pi}}-T\right) + \int_{1-b'}^{1-b} \arg{\zeta(\sigma+iT)}d\sigma + O_{t_0}(\delta).
\end{equation*}

Substituting this into \eqref{eq:lit'}, we have
\begin{equation} \label{eq:n2t}
N_2(T) = \frac{T}{2\pi}\log{\frac{T}{4\pi}} - \frac{T}{2\pi} + \frac{1}{\pi\delta}\left\{\int_b^{b'} \arg{G_2(\sigma+iT)}d\sigma + \int_{1-b'}^{1-b} \arg{\zeta(\sigma+iT)}d\sigma \right\}+ O_{t_0}(1).
\end{equation}

Taking the limit $\delta\rightarrow0$, by the mean value theorem,
\begin{equation*}
\lim_{\delta\rightarrow0} \frac{1}{\pi\delta}\int_b^{b'} \arg{G_2(\sigma+iT)}d\sigma = \frac{1}{2\pi}\arg{G_2\left(\frac{1}{2}+iT\right)}
\end{equation*}
by noting that $b=\frac{1}{2}-\delta$ and $b'=\frac{1}{2}-\frac{\delta}{2}$.
And similarly,
\begin{equation*}
\lim_{\delta\rightarrow0} \frac{1}{\pi\delta}\int_{1-b'}^{1-b} \arg{\zeta(\sigma+iT)}d\sigma = \frac{1}{2\pi}\arg{\zeta\left(\frac{1}{2}+iT\right)}.
\end{equation*}
Substituting these into equation \eqref{eq:n2t}, we have
\begin{equation*}
N_2(T) = \frac{T}{2\pi}\log{\frac{T}{4\pi}} - \frac{T}{2\pi} + \frac{1}{2\pi}\arg{G_2\left(\frac{1}{2}+iT\right)} + \frac{1}{2\pi}\arg{\zeta\left(\frac{1}{2}+iT\right)} + O_{t_0}(1).
\end{equation*}

If $2\leq T<t_0$, then $N_2(T)\leq N_2(t_0) = O_{t_0}(1)$.
Hence the above equation holds for any $T\geq2$ which satisfies the conditions of Proposition \ref{prop6}.\\
\end{proof}

\vspace{3mm}
\noindent
\underline{\textbf{Proof of Theorem \ref{cha2}:}}\\

Firstly we consider for $T\geq2$ which satisfies $\zeta''(\sigma+iT)\neq0$ and $\zeta(\sigma+iT)\neq0$ for any $\sigma\in\mathbb{R}$.
By Lemma \ref{lem5},
\begin{equation*}
\arg{G_2\left(\frac{1}{2}+iT\right)} = O_a\left(\frac{\log T}{(\log{\log{T}})^{1/2}}\right)
\end{equation*}
and again from equation (2.23) of \cite[p. 2251]{aka}, we have
\begin{equation*}
\arg{\zeta\left(\frac{1}{2}+iT\right)} = O\left(\frac{\log T}{\log{\log{T}}}\right).
\end{equation*}
Substituting these into Proposition \ref{prop6}, we obtain
\begin{equation*}
N_2(T) = \frac{T}{2\pi}\log{\frac{T}{4\pi}} - \frac{T}{2\pi} + O_{a,t_0}\left(\frac{\log{T}}{(\log{\log{T}})^{1/2}}\right).
\end{equation*}

Next, if $\zeta(\sigma+iT)=0$ or $\zeta''(\sigma+iT)=0$ for some $\sigma\in\mathbb{R} \quad (T\geq2)$, then we again take a small $0<\epsilon<1$ such that $\zeta''(\sigma+i(T\pm\epsilon))\neq0$ and $\zeta(\sigma+i(T\pm\epsilon))\neq0$ for any $\sigma\in\mathbb{R}$ as in the proof of Theorem \ref{cha1}.
Then similarly, we can show that the increment of the value of $N_2(T)$ can be included in the error term of the above equation.\\

Finally, since $a$ and $t_0$ are fixed constants,
\begin{equation*}
N_2(T) = \frac{T}{2\pi}\log{\frac{T}{4\pi}} - \frac{T}{2\pi} + O\left(\frac{\log{T}}{(\log{\log{T}})^{1/2}}\right)
\end{equation*}
holds for any $T\geq2$.\\
\qed


\section{Final Remarks}
\label{sec:4}

In this paper we have proven that we can extend the results of \cite{aka} to the second derivative of the Riemann zeta function.
However, it is favorable to extend these results to the general case, that is to the $k$-th derivative of the Riemann zeta function (for $k\geq1$, $k\in\mathbb{Z}$) which is our main goal at present.

Finally, the author would like to dedicate special thanks to Prof. Kohji Matsumoto and senior Ryo Tanaka for their valuable advices along the way.


\end{document}